# Séries Gevrey de type arithmétique, II. Transcendance sans transcendance

Par Yves André

### Introduction

Dans ce second volet, nous étudions les propriétés diophantiennes des valeurs de séries Gevrey de type arithmétique d'ordre non nul, en des points algébriques. Nous nous fondons sur le fait, prouvé dans le premier volet, que *l'opérateur différentiel d'ordre minimal annulant une telle série n'a pas de singularités non triviales en dehors de l'origine et de l'infini.*

Nous montrons comment tirer de ce fait des propriétés de transcendance, et retrouver en particulier le théorème fondamental de la théorie de Siegel-Shidlovskii sur l'indépendance algébrique des valeurs de $E$-fonctions en des points algébriques.

Le paradoxe du titre marque le fort constraste entre l'aspect *qualitatif* de ce nouvel argument, et le caractère éminemment quantitatif des techniques de transcendance traditionnelles (ajustement de paramètres auxiliaires...).

Dans le cas particulier des valeurs de la fonction exponentielle, nous obtenons une (ou plutôt deux) preuve du théorème de Lindemann-Weierstrass proche de celle "$p$-adique" de Bézivin-Robba (et qui en est inspirée), où la transformation de Laplace joue un rôle-clé.

Nous traitons de manière analogue le cas des séries Gevrey d'ordre positif, où l'on rencontre d'"étranges" formules comme $\sum_{n\geq 0} n \cdot n! = -1$, valable $p$-adiquement pour tout $p$.

Enfin, nous discutons brièvement le problème des $q$-analogues, et prouvons à titre d'illustration deux $q$-analogues du théorème de Lindemann-Weierstrass (concernant les valeurs de fonctions thêta et de la $q$-exponentielle respectivement).

### Plan
1. Sur le théorème de Lindemann-Weierstrass
2. Valeurs de séries Gevrey d'ordre $< 0$ de type arithmétique (e.g. $E$-fonctions)
3. Valeurs de séries Gevrey d'ordre $> 0$ de type arithmétique (e.g. 3-fonctions)
4. Deux $q$-analogues du théorème de Lindemann-Weierstrass

N.B. Les références au premier volet de cet article seront précédées de I.

*Je remercie D. Bertrand et J. P. Ramis de leur intérêt et pour d'utiles commentaires.*

# 1. Sur le théorème de Lindemann-Weierstrass

Il s'agit des énoncés équivalents suivants, qui contiennent la transcendance de $e$ et de $\pi$:

i) *Si $a_1, \ldots, a_m$ sont des nombres algébriques distincts, les nombres complexes $e^{a_1}, \ldots, e^{a_m}$ sont linéairement indépendants sur $\overline{\mathbb{Q}}$ (la clôture algébrique de $\mathbb{Q}$ dans $\mathbb{C}$);*

ii) *Si $a_1, \ldots, a_m$ sont des nombres algébriques linéairement indépendants sur $\overline{\mathbb{Q}}$, les nombres complexes $e^{a_1}, \ldots, e^{a_m}$ sont algébriquement indépendants sur $\mathbb{Q}$.*

Si $K$ est la clôture galoisienne du corps de nombres engendré par $a_1, \ldots, a_m$ et par des nombres algébriques $b_1, \ldots, b_m$, on remarque que la nullité de $\Sigma b_i e^{a_i}$ entraîne l'annulation au point 1 du polynôme exponentiel

$$\prod_{\sigma \in \mathrm{Gal}(K/\mathbb{Q})} (\Sigma b_i^\sigma e^{a_i^\sigma z}) \in \mathbb{Q}[[z]].$$

On en déduit que les énoncés précédents équivalent encore à chacun des suivants:

iii) *Soit $F \in \mathbb{Q}[[z]]$ un polynôme exponentiel qui s'annule en 1; alors $G = (z-1)^{-1} F$ est un polynôme exponentiel.*

iv) *Soient $a_1, \ldots, a_m$ sont des nombres algébriques distincts, et $b_1, \ldots, b_m$ des nombres algébriques tels que $\Sigma b_i e^{a_i z} \in \mathbb{Q}[[z]]$; alors si $\Sigma b_i e^{a_i} = 0$, on a $b_1 = \cdots = b_m = 0$.*

Prouvons iii). Remarquons d'abord que $G$, tout comme $F$, est le développement à l'origine d'une fonction entière d'ordre exponentiel 1, et est *holonome* (i.e. est solution d'une équation différentielle linéaire à coefficients polynômiaux). De plus, les coefficients de $F = \sum_{n \geq 0} \frac{a_n}{n!} z^n$ et de $G = \sum_{n \geq 0} \frac{b_n}{n!} z^n$ sont liés par la formule $b_n = -\sum_{k=0}^{n} \frac{n!}{k!} a_k$ qui montre que la suite des $b_n$ satisfait, tout comme la suite des $a_n$, à la condition $(\mathbf{G})_{\mathrm{dén}}$ (cf. I.1.1). On en conclut que $G$ est une $E$-fonction. Considérons alors les $G$-fonctions $f$ et $g$ associées à $F$ et $G$ respectivement. Puisque $F$ est un polynôme exponentiel, $f = \sum_{n \geq 0} a_n z^n$ est une fonction rationnelle: $f \in \mathbb{Q}(z)$.

Les formules de I.1.4. se traduisent en le tableau de correspondances suivant:

(1.1) $$f \leftrightarrow F,$$

(1.2) $$z \frac{d}{dz} f \leftrightarrow z \frac{d}{dz} F,$$

(1.3) $$\left(z^2 \frac{d}{dz} + z\right) f \leftrightarrow zF.$$

Ainsi, l'équation $F = (z-1)G$ se traduit par l'équation différentielle

$$(1.4) \qquad f = (z^2 \frac{d}{dz} + z - 1)g,$$

soit encore par le système différentiel

$$(1.5) \qquad z^2 \frac{d}{dz} \begin{pmatrix} g \\ 1 \end{pmatrix} = \begin{pmatrix} 1-z & f \\ 0 & 0 \end{pmatrix} \begin{pmatrix} g \\ 1 \end{pmatrix}.$$

Jusqu'à ce point, l'argument suit [BéR]. Remarquons maintenant que 0 est une singularité *irrégulière* de (1.5). Or les théorèmes I.3.2, 3.3, 3.4 (plus précisément, leurs analogues pour les systèmes différentiels [A1]) entraînent que *tout système différentiel à coefficients rationnels, admettant une solution dont les composantes sont des G-fonctions linéairement indépendantes sur $\mathbb{Q}(z)$, est nécessairement fuchsien.* Comme $g$ est une $G$-fonction, il suit de là que $g$ est en fait une fonction rationnelle (cet argument remplace le critère de rationalité de [BéR]). Donc $G$ est un polynôme exponentiel, ce qu'il fallait démontrer.

## 2. Valeurs de séries Gevrey d'ordre $< 0$ de type arithmétique (e.g. $E$-fonctions)

2.1. Soient $K$ un corps de nombres, et $\xi$ un élément de $K$. Soit $F \in K[[z]]$ une série Gevrey d'ordre $s \in \mathbb{Q}_{<0}$ de type arithmétique holonome. Rappelons que cela signifie que $F$ s'écrit $\sum_{n \geq 0} n!^s a_n z^n$, où les nombres algébriques $a_0, a_1, a_2, \ldots$ vérifient la condition

(**G**): Il existe une constante $C > 0$ telle que pour tout $n$

  (**G**)$_{\text{conj}}$: Les conjugués de $a_n$ sont de module inférieur à $C^n$,

  (**G**)$_{\text{dén}}$: Le dénominateur commun à $a_0, a_1, a_2, \ldots a_n$ est inférieur à $C^n$.

Soit $\Theta \in K[z, \frac{d}{dz}]$ un opérateur d'ordre minimal (en $\frac{d}{dz}$) annulant $F$. Voici une première généralisation du théorème de Lindemann-Weierstrass (sous la forme iii) ci-dessus).

THÉORÈME 2.1.1. *Supposons que pour toute place infinie $v$ de $K$, la fonction analytique définie par $F$ s'annule en $\xi$. Alors $\Theta$ admet une base de solutions dans $(z - \xi)K[[z - \xi]]$; en particulier, $\xi$ est une singularité (triviale*[1]*) de $\Theta$ si $F \neq 0$.*

*Preuve.* Prouvons d'abord le

---

[1] Rappelons qu'une singularité $\zeta$ est dite triviale si l'opérateur différentiel admet une base de solutions séries formelles en $z - \zeta$.

LEMME 2.1.2. $G = (z - \xi)^{-1}F$ est une série Gevrey d'ordre $s$ de type arithmétique et holonome.

L'holonomie est claire. Les coefficients de $F = \sum_{n \geq 0} n!^s a_n z^n$ et de $G = \sum_{n \geq 0} n!^s b_n z^n$ sont liés par la formule

$$b_n = -\sum_{k=0}^{n} (\frac{n!}{k!})^{|s|} \xi^{k-n-1} a_k$$

qui montre que la suite des $b_n$ satisfait, tout comme la suite des $a_n$, à la condition $(\mathbf{G})_{\text{dén}}$. D'autre part, l'hypothèse d'annulation entraîne que pour toute plongement complexe de $K$, la fonction analytique définie par $G$ est (tout comme celle définie par $F$) entière et d'ordre exponentiel $|\frac{1}{s}|$. On en conclut que $G$ est une série Gevrey d'ordre $s$ de type arithmétique (cf. I.1.2, remarque c)), ce qui démontre le lemme.

Le point-clé de la preuve de 2.1.1. est alors qu'*un opérateur d'ordre minimal annulant une série Gevrey d'ordre $s \neq 0$ de type arithmétique n'a pas de singularité non-triviale en dehors de $0$ et de l'infini* (cf. I. introduction, et plus spécifiquement I.4.4 dans le cas des $E$-fonctions). Pour $G$, un tel opérateur peut être pris égal au produit $\Theta \cdot (z - \xi)$. On en déduit que $\Theta$ admet une base de solutions dans $(z - \xi) K[[z - \xi]]$.  C.Q.F.D

Par application itérée de 2.1.1, et en faisant $s = -1$, on obtient le résultat suivant:

COROLLAIRE 2.1.3. *Soit $F \in \mathbb{Q}[[z]]$ une $E$-fonction s'annulant à l'ordre $N > 0$ en $1$ et soit $\Theta \in \mathbb{Q}[z, \frac{d}{dz}]$ un opérateur d'ordre minimal annulant $F$. Alors $\Theta$ admet une base de solutions dans $(z-1)^N K[[z-1]]$.*

2.2. *Application: une autre preuve de Lindemann-Weierstrass.* Déduisons la variante 1.iv) de Lindemann-Weierstrass du théorème 2.1.1. On prend pour $E$-fonction $F = \Sigma b_i e^{a_i z} \in \mathbb{Q}[[z]]$ ($K = \mathbb{Q}$). On peut donc choisir pour $\Theta$ un opérateur différentiel à coefficients constants. Si $F$ s'annule en $\xi = 1$, (2.1.1) entraîne que $\Theta$ a une singularité en $1$ si $F$ n'est pas identiquement nulle. C'est impossible, donc $F = b_1 = \cdots = b_m = 0$.

2.3. *Application: le théorème de Siegel-Shidlovskii.* Montrons comment 2.1.3 permet de retrouver le théorème fondamental de la théorie de Siegel-Shidlovskii ([Si], [Sh, p. 139]):

THÉORÈME 2.3.1. *Soient $E_1, \ldots, E_\mu$ des $E$-fonctions vérifiant un système différentiel $\frac{d}{dz} E_i = \sum_{j=1}^{\mu} A_{i,j} E_j$, où $A_{i,j} \in \overline{\mathbb{Q}}(z)$. Soit $\xi$ un nombre algébrique non nul, distinct des pôles des $A_{i,j}$. Alors le degré de transcendance homogène de $E_1(\xi), \ldots, E_\mu(\xi)$ sur $\overline{\mathbb{Q}}$ est égal au degré de transcendance homogène de $E_1, \ldots, E_m$ sur $\overline{\mathbb{Q}}(z)$.*

Nous renvoyons à [Sh] pour la déclinaison des nombreux énoncés de transcendance et d'indépendance algébrique de valeurs de fonctions classiques qu'on obtient en spécialisant cet énoncé.

*Preuve.* Quitte à remplacer la variable $z$ par $\xi z$, on peut supposer que $\xi = 1$. Il est clair que le degré de transcendance homogène $d_1$ de $E_1(1), \ldots, E_\mu(1)$ sur $\overline{\mathbb{Q}}$ est au plus égal au degré de transcendance homogène $d$ de $E_1, \ldots, E_\mu$ sur $\overline{\mathbb{Q}}(z)$. Soit $K$ un corps de nombres galoisien sur $\mathbb{Q}$ contenant les coefficients des $E_i$, tel que $A_{i,j} \in K(z)$ et que l'idéal des relations homogènes entre les $E_i(1)$ admette des générateurs dans $K[T_1, \ldots, T_\mu]$. Supposons, par l'absurde, que $d_1 < d$.

Il existe alors une relation homogène $Q(E_1(1), \ldots, E_\mu(1)) = 0$ de degré $\delta$ "fortement non-triviale" au sens de [Bo, 12]. Fortement non-triviale veut dire ceci: soient $K[z]_1$ la localisation en 1 de $K[z]$, $I \subset K[z]_1[T_1, \ldots, T_\mu]$ l'idéal homogène des relations homogènes entre les $E_i$, $X = \operatorname{Proj} K[z]_1[T_1, \ldots, T_\mu]/I$, et $X_1$ la fibre spéciale de $X$; alors $Q$ définit une hypersurface de $X_1$ (en d'autres termes, $\dim(X_1 \cap \{Q = 0\}) < d$). Bien entendu, on suppose que $\mathcal{O}_X$ est sans $K[z]_1$-torsion (i.e. on divise les relations homogènes par $z - 1$ autant que faire se peut); ainsi $X$ est un schéma intègre, mais sa fibre spéciale $X_1$ peut n'être ni irréductible ni réduite.

Soit $k$ un entier naturel. Comme $\mathcal{O}_X$ est sans $K[z]_1$-torsion, le $K[z]_1$-module de type fini $H^\circ(X, \mathcal{O}(k))$ est libre; son rang est donné par la fonction de Hilbert $m = h(X_{K(z)}, k) = h(X_1, k)$. Pour $k$ assez grand, l'application

$$H^\circ(\operatorname{Proj} K[z]_1[T_1, \ldots, T_\mu], \mathcal{O}(k)) \to H^\circ(X, \mathcal{O}(k))$$

est surjective. Vu comme module des polynômes de degré $k$ en les $E_i$, $H^\circ(X, \mathcal{O}(k))$ est ainsi naturellement muni d'une action de $\frac{d}{dz}$. Les éléments $F_1, \ldots, F_m$ d'une base fixée de $H^\circ(X, \mathcal{O}(k))$ (formée de polynômes de degré $k$ en les $E_i$ linéairement indépendants sur $K(z)$) vérifient donc un système différentiel

$$(*) \qquad \frac{d}{dz} F_j = \sum_{h=1}^{m} B_{j,h} F_h, \quad \text{où } B_{j,h} \in K[z]_1.$$

La dimension $n$ du $K$-espace $H^\circ(X_1 \cap \{Q = 0\}, \mathcal{O}(k))$ engendré par les $F_j(1)$ est donnée par la fonction de Hilbert $n = h(X_1 \cap \{Q = 0\}, k)$. Pour $k \to \infty$, $n/m = h(X_1 \cap \{Q = 0\}, k)/h(X_1, k)$ est équivalent à $\delta d/k$. En particulier, on a pour $k$ assez grand: $n/m = (1 - \varepsilon)/[K : \mathbb{Q}]$, avec $\varepsilon > 0$. Dans la suite, un tel entier $k$ (et donc $\varepsilon$) sera fixé.

Soient $Z_1 = (Z_{1,1}, \ldots, Z_{1,m}), \ldots, Z_n$ des solutions de $(*)$ dans $K[[z-1]]^m$ dont les valeurs en 1 forment une base de $H^\circ(X_1 \cap \{Q = 0\}, \mathcal{O}(k))$, et soit $N$ un entier naturel. Il existe des polynômes $P_1, \ldots, P_m \in K[z]$ non tous nuls, tels que

i) $\mathrm{ord}_1 \sum_{h=1}^m P_h Z_{j,h} \geq N$ pour $j = 1, \ldots, n$,

ii) $\deg P_h \leq nN/m = (1-\varepsilon)N/[K:\mathbb{Q}]$ pour $h = 1, \ldots m$.

En effet, la condition i) s'exprime comme un système linéaire à $nN$ équations en $m(1 + \max \deg P_h)$ inconnues, les coefficients des $P_h$. Compte tenu de ce que 1 n'est pas une singularité de $(*)$, on déduit de i) que $\mathrm{ord}_1 \sum_{h=1}^m P_h F_h \geq N$ (ici $\sum_{h=1}^m P_h F_h$ est considérée comme fonction entière via le plongement fixé de $K$ dans $\mathbb{C}$).

Considérons la série

$$F = \prod_{\sigma \in \mathrm{Gal}(K/\mathbb{Q})} \left( \sum_{h=1}^m P_h^\sigma F_h^\sigma \right) \in \mathbb{Q}[[z]].$$

C'est une $E$-fonction qui s'annule en 1 à l'ordre $N$. Si $\Theta \in \mathbb{Q}[z, \frac{d}{dz}]$ est un opérateur d'ordre minimal annulant $F$, alors d'après 2.1.3, les exposants de $\Theta$ en 1 sont tous des entiers $\geq N$. D'autre part, $F$ s'écrit $\sum_{i=1}^{m^{[K:\mathbb{Q}]}} Q_i Y_i$, où les $Q_i \in K[z]$ sont de degré $\leq (1-\varepsilon)N$ et où les $Y_i$ (i.e. les monômes $\prod_k F_k^{\sigma_k}$) sont les composantes d'une solution, indépendante de $N$, d'un système différentiel à coefficients dans $K(z)$.

Pour $N$ assez grand, ceci contredit l'inégalité de Fuchs pour les exposants de $\Theta$, établie dans le cas irrégulier dans [BB].[2] On conclut que $d_1 = d$. C.Q.F.D.

*Remarque* 2.4. En fait, le théorème de Siegel-Shidlovskii vaut aussi pour les $E$-fonctions au sens large, i.e. dont les coefficients de Taylor $\frac{a_n}{n!}$ vérifient seulement la condition affaiblie (cf. I.2.3)

($\mathbf{G}^-$): Pour tout $\varepsilon > 0$ et pour tout $n$ assez grand,

   ($\mathbf{G}^-$)$_{\mathrm{conj}}$: Les conjugués de $a_n$ sont de module inférieur à $(n!)^\varepsilon$,

   ($\mathbf{G}^-$)$_{\mathrm{dén}}$: Le dénominateur commun à $a_0, a_1, a_2, \ldots a_n$ est inférieur à $(n!)^\varepsilon$.

La preuve ci-dessus s'étend à ce cas. Le point est qu'*un opérateur d'ordre minimal annulant une $E$-fonction au sens large n'a pas de singularité non-triviale en dehors de $0$ et de l'infini.* Comme on l'a vu en I.5, ceci découle de ce que *la singularités en l'infini d'un opérateur $\phi$ d'ordre minimal annulant une $G$-fonction au sens large est régulière.*

Nous serons très bref sur ce point, puisque la distinction entre les sens larges et stricts est conjecturalement illusoire, et que dans les applications, les $E$-fonctions intervenant le sont clairement au sens strict.

---

[2] Une analyse plus détaillée de ce point se trouve dans le récent survol de D. Bertrand "On André's proof of the Siegel-Shidlovsky theorem", compte-rendu du Colloque franco-japonais de théorie des nombres 1999, Univ. Keio.

Pour établir le point ci-dessus, le critère $p$-adique de régularité de Katz montre qu'il suffit d'établir que $\sum_{p(v)\leq n} \log R_v(\phi,1) = o(\log n)$. Avec les notations de [A, VI], on voit facilement que cette condition découle d'une estimation

$$\sum_{v\,\text{finie}} h_{v,n}(\phi) = \sum_{p(v)\leq n} h_{v,n}(\phi) = o(\log n);$$

or les estimations de loc. cit., p. 122 donnent

$$\sum_{v\,\text{finie}} h_{v,n(\phi)} \leq C_1 \frac{1}{n} \sum_v \log\,\max(1,|a_0|_v,\ldots,|a_{nC_2}|_v) + C_3$$

(pour des constantes $C_i$ indépendantes de $n$), et la condition $(\mathbf{G}^-)$ équivaut à

$$\frac{1}{n}\sum_v \log\max(1,|a_0|_{v'},\ldots,|a_n|_v = o(\log n).$$

## 3. Valeurs de séries Gevrey d'ordre $> 0$ de type arithmétique (e.g. $\mathfrak{z}$-fonctions)

3.1. Soient $K$ un corps de nombres et $\phi$ un élément de $K[z,\frac{d}{dz}]$. Rappelons la définition des rayons de solubilité générique $R_v(\phi,1)$: pour toute place finie $v$ de $K$, de caractéristique résiduelle $p = p(v)$, normalisons la valeur absolue $v$-adique par $|p|_v = p^{-1}$; pour tout $r > 0$, $R_v(\phi,1)$ est le rayon de convergence, limité supérieurement à 1 par convention, d'une base de solutions de $\phi$ au point générique. Rappelons encore (I. 3.3) qu'une condition nécessaire et suffisante pour que $\phi$ soit un G-opérateur est:

(3.1.1) $$\prod_{v\,\text{finie}} R_v(\phi,1) \neq 0.$$

Soit d'autre part $g \in K[[z]]$ une $G$-fonction. Alors les rayons de convergence $v$-adiques de $g$, limités supérieurement à 1 par convention, vérifient:

(3.1.2) $$\prod_v R_v(g) \neq 0.$$

See [A, p. 126] (ici, le produit porte sur toutes les places de $K$, finies ou infinies). La non-nullité de $\prod_{v\,\text{finie}} R_v(g)$ se déduit d'ailleurs de I.3.2, 3.3, appliqués à un opérateur d'ordre minimal $\phi$ annulant $g$, et d'un théorème de spécialisation pour les rayons de convergence $v$-adiques; la non-nullité de $\prod_{v\,\text{infinie}} R_v(g)$ se déduit du fait que 0 est une singularité régulière de $\phi$.

La réciproque est l'une des principales questions ouvertes de la théorie des $G$-fonctions:

CONJECTURE 3.1.3.   *Si $g \in K[[z]]$ est holonome et vérifie $\prod_v R_v(g) \neq 0$, alors $g$ est une G-fonction.*

3.2. Soit $\xi$ un élément de $K$. Soit $\mathfrak{f} \in K[[z]]$ une série Gevrey d'ordre $s \in \mathbb{Q}_{>0}$ de type arithmétique, holonome, et soit $\Theta \in K[z, \frac{d}{dz}]$ un opérateur d'ordre minimal annulant $\mathfrak{f}$.

THÉORÈME 3.2.1.  *Supposons la conjecture* 3.1.3 *vraie. Supposons que pour toute place finie $v$ sauf peut-être un nombre fini, la fonction $v$-adique définie par $\mathfrak{f}$ s'annule en $\xi$. Alors $\xi$ est une singularité (triviale) de $\Theta$. En particulier, $\mathfrak{f}(\xi) = 0$ en toute place finie $v$ telle que $\mathfrak{f}(\xi)$ converge, et en toute place à l'infini, la $\frac{1}{s}$-sommation[3] de $\mathfrak{f}$ s'annule en $\xi$.*

*Preuve.* Prouvons d'abord le

LEMME 3.2.2.  $\mathfrak{g} = (z - \xi)^{-1}\mathfrak{f}$ *est une série Gevrey d'ordre $s$ de type arithmétique et holonome.*

L'holonomie est claire. Posons $s = p/q$. Les coefficients de

$$\mathfrak{f} = \sum_{n \geq 0} [\frac{n}{q}]!^p a_n z^n$$

et de

$$\mathfrak{g} = \sum_{n \geq 0} \left[\frac{n}{q}\right]!^p b_n z^n$$

sont liés par la formule

$$b_n = -\sum_{k=0}^{n} \left(\frac{[k/q]!}{[n/q]!}\right)^p \xi^{k-n-1} a_k$$

qui montre que la suite des $b_n$ satisfait, tout comme la suite des $a_n$, à la condition $(\mathbf{G})_{\text{conj}}$ (cf. I.1.1), et que les $R_v(\mathfrak{g})$ sont non nuls pour toute place finie $v$. De plus, l'hypothèse d'annulation montre que $R_v(\mathfrak{g}) = R_v(\mathfrak{f})$ pour presque toute place finie $v$.

D'autre part, considérons la $G$-fonction $f = \sum_{n \geq 0} a_n z^n$ associée à $\mathfrak{f}$ et la série holonome $g = \sum_{n \geq 0} b_n z^n$ associée à $\mathfrak{g}$ (cf. I.2.1). On a alors $R_v(f) = R_v(g)$ pour presque toute place finie $v$, et les $R_v(g)$ sont non nuls pour toute place $v$ de $K$. En particulier $\prod_v R_v(g) \neq 0$. D'après 3.1.3, ceci entraîne que g est une $G$-fonction, donc que $\mathfrak{g}$ est une série Gevrey d'ordre $s$ de type arithmétique.

A partir de là, la preuve de ce que $\xi$ est une singularité triviale de $\Theta$ s'effectue comme en 2.1.1. Ainsi $\Theta$ admet une base de solutions dans $(z - \xi)K[[z - \xi]]$, et 3.2.1 en découle.

---

[3] Si $\overrightarrow{0\xi}$ est une direction singulière, il s'agit de la sommation médiane.

*Exemple* 3.3. On a l'égalité $\sum_{n\geq 0} n \cdot n! = -1$ dans $\mathbb{Z}_p$ pour tout nombre premier $p$ (cf. [Sc, p. 63]). On vérifie alors aisément que la 3-fonction $\mathfrak{f} = 1 + \sum_{n\geq 0} n \cdot n! z^n$ satisfait une équation différentielle inhomogène d'ordre 1, et que $\xi = 1$ en est une singularité triviale. (Dans cet exemple, $\overrightarrow{0\xi}$ est une direction singulière.) Pour d'autres exemples dûs à L. van Hamme et à B. Dragovic, e.g. $= \sum_{n\geq 0} n^5 \cdot (n+1)! = 26$, voir [Sc] et [D].

Ces résultats sont à rapprocher de travaux de Chirskii, où la condition que $\mathfrak{f}$ s'annule en $\xi$ $v$-adiquement pour presque tout $v$ est remplacée par une condition effective d'annulation de $\mathfrak{f}$ en $\xi$ pour tout $v$ de caractéristique résiduelle inférieure à une constante dépendant des données [C].

*Remarques* 3.4. i) La conjecture 3.1.3 entraîne plus généralement une description conjecturale $p$-adique des séries Gevrey holonomes de type arithmétique, qui est à rapprocher de I.4.7. Les détails de cette déduction sont laissés au lecteur.

CONJECTURE 3.4.1. *Soit $s$ un nombre rationnel. Si $f \in K[[z]]$ est holonome, définit une série Gevrey d'ordre $s$ en chaque place à l'infini, et vérifie $\prod_{v \text{ finie}}(R_v(f) \cdot p_v^{\overline{\frac{s}{1-p_v}}}) \neq 0$, alors $f$ est une série Gevrey d'ordre $s$ de type arithmétique.*

(La condition d'holonomie n'est pas superflue, déjà dans le cas $s = 0$; cf. [A, I]).

ii) Il serait intéressant de voir si 3.2.1 subsiste en remplaçant l'hypothèse que *pour toute place finie $v$ sauf peut-être un nombre fini, la fonction $v$-adique définie par $\mathfrak{f}$ s'annule en $\xi$* par: *en toute place à l'infini, la $\frac{1}{s}$-sommation de $\mathfrak{f}$ s'annule en $\xi$*. Le théorème de dualité du volet I est peut-être utile à cet égard.

iii) Il serait intéressant d'étudier avec ces méthodes le cas des $G$-fonctions ($s = 0$). Bien entendu, les arguments qualitatifs utilisés ci-dessus ne suffisent plus, et les propriétés arithmétiques de $\xi$ doivent entrer en ligne de compte (cf. ([A2]). Le cas de séries appartenant à l'algèbre des séries holonomes et Gevrey, d'ordre $s < 0$ variable, de type arithmétique se ramène essentiellement au cas des $E$-fonctions; beaucoup plus difficile, sans doute, serait l'étude du cas où l'on tolère $s \leq 0$ (par exemple la somme d'une $G$-fonction et d'une $E$-fonction).

## 4. Deux $q$-analogues du théorème de Lindemann-Weierstrass

4.1. Soit $q$ un nombre complexe de module $> 1$, et $s$ un nombre rationnel non nul. La notion de série $q$-Gevrey a été introduite par J. P. Bézivin [B]. Une série à coefficients complexes $F = \sum_{n\geq 0} a_n z^n$ est dite *$q$-Gevrey d'ordre $s$* (resp. *d'ordre précis s*) si la série $\sum_{n\geq 0} q^{-sn(n-1)/2}$ a un rayon de convergence non nul (resp. et fini).

Le résultat principal de [Bé] est *que toute série formelle solution d'un opérateur linéaire aux q-différences à coefficients polynômiaux est ou bien convergente, ou bien est q-Gevrey d'ordre précis s*, où $s$ est l'inverse de l'une des pentes (rationnelles) du polygône de Newton. Ce résultat se transpose d'ailleurs tel quel au domaine ultramétrique (noter le contraste avec la notion de série Gevrey, qui dans le cas ultramétrique coïncide avec la notion banale de série convergente).

4.2. Notons $\sigma_q$ la dilatation: $\sigma_q f(z) = f(qz)$, et $\delta_q$ le $q$-analogue de la dérivation: $\delta_q f(z) = \frac{f(qz)-f(z)}{(q-1)z}$. Suivant l'usage, posons $n_q = \frac{1-q^n}{1-q}$, $n_q! = \prod_{m=1}^n m_q$. On a $|n_q!| \approx |q^{n(n-1)/2}|$ quand $n \to \infty$, donc on peut remplacer

$$\sum_{n\geq 0} q^{-sn(n-1)/2} a_n z^n$$

par

$$\sum_{n\geq 0} n_q!^{-s} a_n z^n$$

dans la définition des séries $q$-Gevrey.

On a d'ailleurs deux candidats $q$-analogues de la transformée de Laplace formelle:

$$F^\# = \sum_{n\geq 0} q^{n(n-1)/2} a_n z^{-n-1}, \text{ et } F^+ = \sum_{n\geq 0} n_q! a_n z^{-n-1}.$$

Par exemple, la fonction de Tschakaloff (thêta tronquée)

$$T_q = \sum_{n\geq 0} q^{-n(n-1)/2} z^n$$

vérifie $T_q^\# = \frac{1}{z-1}$, tandis que le $q$-analogue de l'exponentielle $E_q = \sum_{n\geq 0} \frac{z^n}{n_q!}$ vérifie $E_q^+ = \frac{1}{z-1}$. Notons les formules

$$(zF)^\# = \frac{1}{qz}\sigma_{q^{-1}}(F^\#), \quad (zF)^+ = \frac{-1}{q}\delta_{q^{-1}}(F^+).$$

Du point de vue arithmétique qui nous occupe ici, il importe de remarquer que pour $n$ grand, $n_q!$ et $q^{n(n-1)/2}$ sont de natures totalement différentes. Il semble ainsi loisible, pour $q$ algébrique, de considérer tant $T_q$ que $E_q$ comme des $q$-E-fonctions, mais en des sens bien distincts. Plus généralement, pour un nombre algébrique $q$ non nul qui n'est pas une racine de l'unité, cela suggère deux $q$-analogues de la notion de série Gevrey de type arithmétique:

1) Les *séries $q^\#$-Gevrey d'ordre $s$ de type arithmétique* $\sum_{n\geq 0} a_n z^n$, pour lesquelles $\sum_{n\geq 0} q^{-sn(n-1)/2} a_n z^n$ est Gevrey d'ordre 0 de type arithmétique;

(2) Les *séries $q^+$-Gevrey d'ordre $s$ de type arithmétique* $\sum_{n\geq 0} a_n z^n$, pour lesquelles $\sum_{n\geq 0} n_q!^{-s} a_n z^n$ est Gevrey d'ordre 0 de type arithmétique.

Pour toute place $v$ telle que $|q|_v > 1$ (il en existe), une telle série $v$-adique est $q$-Gevrey d'ordre $s$ au sens de Bézivin. Ces définitions sont provisoires; il y aura peut-être lieu de les modifier pour intégrer le cas des séries $q$-hypergéométriques de Heine

$$\begin{aligned}{}_2\Phi_1(a,b,c;q;z) = {}& 1 + \frac{(1-q^a)(1-q^b)}{(1-q^c)}\frac{z}{1-q} \\ & + \frac{(1-q^a)(1-q^{a+1})(1-q^b)(1-q^{b+1})}{(1-q^c)(1-q^{c+1})}\frac{z^2}{(1-q)(1-q^2)} + \cdots\end{aligned}$$

à paramètres $a, b, c$ rationnels.

4.3. La question de savoir si les résultats que nous avons obtenus à propos des séries Gevrey de type arithmétique se transposent aux $q$-analogues est ouverte. Dans ce sens, nous allons établir un analogue du théorème de Lindemann-Weierstrass (variante iv) du §1 ci-dessus) pour les valeurs de $T_q$ et de $E_q$ respectivement (considérées comme $q^\#$- et $q^+$-variantes de l'exponentielle). Nous suivrons pour cela l'argument de Bézivin-Robba.[4]

4.4. Soient $q, \xi, \alpha_1, \ldots \alpha_m, \beta_0, \beta_1 \ldots, \beta_m$ des éléments d'un corps de nombres $K$. On suppose $q, \xi, \alpha_1, \ldots, \alpha_m$ non nuls, et que *q n'est pas une racine de l'unité*. On suppose en outre que *les classes de $\alpha_1, \ldots, \alpha_m$ dans $K^*/q^{\mathbf{Z}}$ soient deux à deux distinctes*.

THÉORÈME 4.4.1. *Si pour toute place $v$ de $K$ telle que $|q|_v > 1$, la fonction $v$-adique $\beta_0 + \beta_q T_q(\alpha_1 z) + \cdots + \beta_m T_q(\alpha_m z)$ s'annule en $\xi$, alors $\beta_0 = \cdots = \beta_m = 0$.*

*Preuve.* Posons $F = \beta_0 + \beta_q T_q(\alpha_1 z) + \cdots + \beta_m T_q(\alpha_m z)$ et $G = (z-\xi)^{-1} F$. Prouvons d'abord le

LEMME 4.4.2. $G^\# \in K\{\frac{1}{z}\}_0^A$, *i.e. est une série Gevrey d'ordre 0 de type arithmétique en la variable $\frac{1}{z}$. De plus, les coefficients de $G^\#$ sont $w$-entiers pour presque toute place finie $w$ de $K$.*

---

[4] Depuis la soumission de cet article, L. di Vizio a démontré bon nombre de $q$-analogues des résultats ci-dessus, notamment un $q$-analogue de 2.1.1 qui permet de retrouver 4.4.1 et 4.5.1 par l'argument de 2.2.

Pour toute place $w$ en dehors d'un ensemble fini, le coefficient $a_n$ de $z^{-n-1}$ dans $F^\#$ est $w$-entier puisque $F^\#$ est une fonction rationnelle

$$\left(F^\# = \frac{\beta_0}{z} + \frac{\beta_1}{z - \alpha_1} + \cdots + \frac{\beta_m}{z - \alpha_m}\right).$$

Le coefficient de $z^{-n-2}$ dans $G^\#$ s'écrit

$$\sum_{k=0}^{n} \sqrt{q}^{(n(n-1)-k(k-1))} \xi^{n-k} \left(\sum_{j=0}^{m} \beta_j \alpha_j^k\right),$$

où par convention on a posé $\alpha_0^k = 1$ si $k = 0$, $\alpha_0^k = 0$ sinon. Ce coefficient est donc $w$-entier pour tout $w$ hors d'un ensemble fini indépendant de $n$. On voit aussi que la série $G^\#$ converge $w$-adiquement pour toute place $w$ telle que $|q|_w \leq 1$. Il ne reste plus qu'à établir que $G^\#$ converge $v$-adiquement aussi pour toute place $v$ telle que $|q|_v > 1$. Or pour ces places, on a par hypothèse l'égalité de fonctions $v$-adiques $F = (z-\xi)G$, d'où $F^\# = (\frac{1}{qz}\sigma_q^{-1} - \xi)(G^\#)$. Ceci fournit une équation aux $q$-différences du second ordre à coefficients polynômiaux pour $G^\#$, dont les pentes à l'infini sont 0 et 1. On déduit alors du théorème de Bézivin (4.1) que si $G^\#$ diverge $v$-adiquement, alors $G^\#$ est $q$-Gevrey d'ordre précis 1. Ainsi $G$ serait Gevrey d'ordre précis 0. C'est impossible car $G$ est, tout comme $F$, une fonction entière. Donc $G^\#$ converge $v$-adiquement, et ceci achève la preuve du lemme.

LEMME 4.4.3. *Soit $w$ une place pour laquelle $|q|_w < 1$ (il en existe d'après Kronecker). Alors $G^\#$ est méromorphe sur tout le "plan" $w$-adique privé de l'origine.*

Fixons un entier $n_0$. On peut décomposer $G^\#$ en deux termes:

$$G^\# = P_{n_0}(\sigma_{q^{-1}})F^\# + \sum_{n \geq n_0} \sum_{k=0}^{n-n_0} \sqrt{q}^{(n(n-1)-k(k-1))} \xi^{n-k} \left(\sum_{j=0}^{m} \beta_j \alpha_j^k\right) z^{-n-2}.$$

Le premier terme est une fonction rationnelle, tandis que le second est une série convergente pour

$$|z|_w > |q|_w^{n_0} \cdot \mathrm{Max}\left(|\alpha_1|_w, \ldots, |\alpha_m|_w, |\xi|_w\right).$$

On en déduit le lemme en prenant $n_0$ arbitrairement grand.

Ces lemmes permettent d'appliquer à $G^\#$ le critère de Borel-Dwork [Am, 5.3], qui montre que $G^\#$ est rationnelle. De ce que, par hypothèse, les pôles de $F^\#$ (c'est-à-dire les $\alpha_i$) sont simples et sont dans des classes distinctes dans $K^*/q^{\mathbf{Z}}$, on voit que l'égalité $F^\# = (\frac{1}{qz}\sigma_{q^{-1}} - \xi)(G^\#)$ n'est possible que si $F^\# = 0$, c'est-à-dire $\beta_0 = \cdots = \beta_m = 0$. \hfill C.Q.F.D.

Le théorème 4.4.1 fournit des renseignements sur la nature diophantienne des valeurs de la fonction thêta

$$\theta(z,q) := \sum_{n=-\infty}^{+\infty} q^{n^2} z^n \left( = \theta_3\left(\frac{\log z}{2i}, q\right)\right)$$

dans les notations de [WW]).

COROLLAIRE 4.4.4.  *Supposons que $q$ soit l'inverse d'un entier naturel $> 1$ (resp. un nombre rationnel positif $< 1$ dont le numérateur est une puissance d'un nombre premier $p$). Soient $\xi_1, \ldots, \xi_m, \gamma_1, \ldots, \gamma_m$ des nombres rationnels non nuls. On suppose que les $\xi_m^{\pm 1}, \ldots \xi_m^{\pm 1}$ sont dans des classes distinctes dans $\mathbb{Q}^*/q^{\mathbf{Z}}$. Alors le nombre réel (resp. $p$-adique) $\gamma_1\theta(\xi_1, q) + \cdots + \gamma_m\theta(\xi_{m,q})$ est irrationnel.*

(En effet,

$$\gamma_1\theta(\xi_1, q) + \cdots + \gamma_m\theta(\xi_m, q) = \sum_{i=1}^m \gamma_i \left(T_{q^{-1}}(q\xi_i) + T_{q^{-1}}\left(\frac{q}{\xi_i}\right) - 1\right),$$

et il n'y a qu'une place $v$ à considérer.)

*Remarque* 4.4.5. Il est intéressant de comparer ces résultats avec ceux de [T] et [Bun], obtenus par des méthodes complètement différentes des nôtres. Ces auteurs s'intéressent au cas où $q$ et les arguments $\xi_i$ sont rationnels (ou éventuellement dans un corps quadratique imaginaire), et font l'hypothèse archimédienne que le numérateur de $q$ est beaucoup plus grand que son dénominateur. Dans cette situation particulière, nos résultats ne couvrent les leurs que lorsque $q$ est entier.

4.5. Passons à la $q$-exponentielle, et remarquant d'abord le développement en produit $E_q(z) = \prod_{m=1}^{\infty}[1 + (q-1)\frac{z}{q^m}]$, qui donne les zéros de $E_q$ : $z = \frac{q^m}{1-q}$ pour m entier $> 0$.

Soient de nouveau $q, \xi, \alpha_1, \ldots, \alpha_m, \beta_0, \eta_1, \ldots, \beta_m$ des éléments d'un corps de nombres $K$. On suppose $q, \xi, \alpha_1 \ldots, \alpha_m$ non nuls, et que *pour aucune place archimédienne, $q$ n'est sur le cercle unité.*[5] On suppose en outre que *les classes de $\alpha_1, \ldots, \alpha_m$ dans $K^*/q^{\mathbf{Z}}$ soient deux à deux distinctes.*

THÉORÈME 4.5.1.  *Si pour toute place $v$ de $K$ telle que $|q|_v > 1$, la fonction $v$-adique $\beta_0 + \beta_1 E_q(\alpha_1 z) + \cdots + \beta_m E_q(\alpha_m z)$ s'annule en $\xi$, alors $\beta_0 = \beta_1 E_q(\alpha_1 \xi) = \cdots = \beta_m E_q(\alpha_m \xi) = 0$.*

---

[5] Cette hypothèse s'avère superflue, si l'on tient compte de ce que pour toute place archimédienne $v|q|_v > 1$ telle que le nombre algébrique $q$ soit de valeur absolue 1, le réel $\log |q|_v/i\pi$ n'est pas de Liouville. Ceci entraîne la convergence $v$-adique de la série $G^+$ introduite plus loin.

*Preuve.* Posons $F = \beta_0 + \beta_1 E_q(\alpha_1 z) + \cdots + \beta_m E_q(\alpha_m z)$ et $G = (z-\xi)^{-1}F$. Un argument parallèle à celui de 4.4.2 établit le

LEMME 4.5.2.   $G^+ \in K\{\frac{1}{z}\}_0^A$. *De plus, les coefficients de $G^+$ sont $w$-entiers pour presque toute place finie $w$ de $K$.*

Le seul écart d'avec 4.4.2 est que, le coefficient de $z^{-n-2}$ dans $G^+$ s'écrivant maintenant

$$\sum_{k=0}^n \frac{n_q!}{k_q!}\xi^{n-k}\left(\sum_{j=0}^m \beta_j \alpha_j^k\right),$$

c'est en tenant compte de ce que pour aucune place archimédienne, $q$ n'est sur le cercle unité que l'on voit que la série $G^+$ converge $w$-adiquement pour toute place $w$ telle que $|q|_w \leq 1$.

LEMME 4.5.3.   *Pour presque toute place finie $w$ de $K$, $G^+$ définit une fonction méromorphe $w$-adique sur le domaine $|z|_w > p_w^{\frac{-1}{p_w 1}}$.*

On a

$$F^+ = \left(\frac{-1}{q}\delta_{q^{-1}} - \xi\right)G^+ = \frac{\beta_0}{z-\alpha_0} + \frac{\beta_1}{z-\alpha_1} + \cdots + \frac{\beta_m}{z-\alpha_m}$$

(avec $\alpha_0 = 0$). Suivant l'idée de [BéR], écrivons formellement

$$G^+ = -\xi\left(1 + \frac{q}{\xi q}\delta_{q^{-1}}\right)^{-1}F^+ = \sum_{n\geq 0}(-\xi)^{-n-1}\left(\frac{1}{q}\delta_{q^{-1}}\right)^n\left(\sum_{j=0}^m \frac{\beta_j}{z-\alpha_j}\right).$$

Or par récurrence, on déduit aisément la formule

$$\left(\frac{1}{q}\delta_{q^{-1}}\right)^n \frac{1}{z-\alpha} = \frac{(-1)^n n_q!}{(z-\alpha)\cdots(z-\alpha q^n)}.$$

Comme pour tout premier $p$ et tout multiple m de $p-1$, $p$ divise $q^m - 1$, on voit que $p^{[n/p-1]}$ divise $\prod_{m=1}^n q^m - 1$; donc $\text{val}_w n_q! \geq [n/p_w - 1]$ pour toute place finie $w$ telle que $|q-1|_w \geq 1$. Ceci démontre que pour toute place finie $w$ telle que $q, q-1, \xi$, et les $\alpha_j$ soient des $w$-unités, la série

$$-\sum_{n\geq 0}(\xi)^{-n-1}\left(\sum_{j=0}^m \frac{n_q! \beta_j}{(z-\alpha_j)\cdots(z-\alpha_j q^n)}\right)$$

converge au voisinage de l'infini, est égale à $G^+$, et définit en fait une fonction méromorphe sur le domaine $|z|_w > p_w^{-1/(p_w-1)}$.

Ces lemmes permettent d'appliquer à $G^+$ le critère de Polya-Bertrandias [Am, 5.4.6], qui montre que $G^+$ est rationnelle, puisque $\prod_w p_w^{-1/(p_w-1)} = 0$.

Comme
$$F^+ = \frac{\beta_0}{z} + \frac{\beta_1}{z - \alpha_1} + \cdots + \frac{\beta_m}{z - \alpha_m} = \left(\frac{-1}{q}\delta_{q^{-1}} - \xi\right) G^+$$
est à pôles simples, il en est de même de $G^+$; de plus les pôles de $G^+$ doivent être de la forme $a_i q^n$, avec $n$ entier. Compte tenu de ce que, par hypothèse, les pôles de $F^+$ sont dans des classes distinctes dans $K^*/q^{\mathbf{Z}}$, on conclut par le lemme suivant:

LEMME 4.5.4. *Pour $\alpha \neq 0$, l'équation*
$$\frac{1}{z - \alpha} = \left(\frac{-1}{q}\delta_{q^{-1}} - \xi\right) \left(\sum_{n=-M}^{n=N} \frac{\gamma_n}{z - \alpha q^n}\right)$$
*n'a lieu que si $\alpha\xi = \frac{q^m}{1-q}$ pour un entier $m > 0$.*

En effet,
$$\left(\frac{-1}{q}\delta_{q^{-1}} - \xi\right) \left(\sum_{n=-M}^{n=N} \frac{\gamma_n}{z - \alpha q^n}\right)$$
$$= \sum_{n=-M}^{n=N} \frac{\gamma_n}{\alpha q^n(1-q)} \left(\frac{1}{z - \alpha q^n} - \frac{1}{z - \alpha q^{n+1}}\right) - \frac{\xi\gamma_n}{z - \alpha q^n}$$
$$= \sum_{n=-M}^{n=N+1} \left[\gamma_m \left(\frac{1}{\alpha q^n(1-q)} - \xi\right) - \frac{\gamma_{m-1}}{\alpha q^n(1-q)}\right] \frac{1}{z - \alpha q^n},$$
en posant $\gamma_{N+1} = \gamma_{-M-1} = 0$. On en tire les équations
$$\gamma_{m-1} = [1 - \alpha\xi q^m(1-q)]\gamma_m \quad \text{pour } m \neq 0,$$
$$\gamma_1 - \gamma_0[1 - \alpha\xi(1-q)] = \alpha(q-1) \neq 0,$$
d'où découle immédiatement le résultat. C.Q.F.D.

*Remarque* 4.5.5. Il est amusant de constater que pour $q =$ une puissance négative d'un nombre premier $p$, on obtient ainsi une variante $p$-adique de l'analogue de Lindemann-Weierstrass pour la $q$-exponentielle, alors qu'aucune version $p$-adique de Lindemann-Weierstrass n'est connue pour l'exponentielle classique.

COROLLAIRE 4.5.6 (R. Wallisser [W]). *Supposons que $q$ soit un entier naturel $> 1$. Soit $\xi$ un nombre rationnel non nul. Alors le nombre réel $E_q(\xi)$ est irrationnel ou nul. Pour $|\xi| < 1$, le nombre réel $E_{q^{-1}}(\xi)$ est irrationnel.*

La première assertion est conséquence immédiate de 4.5.1 avec $v = \infty$. La seconde découle de là, de l'égalité $E_{q^{-1}}E_q(-\xi) = 1$, valable pour $|\xi| < 1$, et de ce que $E_q$ ne s'annule pas dans le disque unité. De même, on obtient:

COROLLAIRE 4.5.7. *Supposons que $q$ soit un nombre rationnel positif $< 1$ dont le numérateur est une puissance d'un nombre premier $p$. Soit $\xi$ un nombre rationnel non nul. Alors le nombre $p$-adique $E_q(\xi)$ est irrationnel ou nul. Pour $|\xi| < 1$, le nombre $p$-adique $E_{q^{-1}}(\xi)$ est irrationnel.*

*Remarque* 4.5.8. La notion de série $q$-Gevrey d'ordre $s$ (et *a fortiori* celle de série $q$-Gevrey de type arithmétique) n'est pas stable par multiplication. Ceci rend difficilement accessible les questions de transcendance des valeurs spéciales de fonctions du type $T_q$ ou $E_q$. Il est toutefois conjecturé que pour $q \in \overline{\mathbb{Q}}$, les valeurs de ces séries en tout nombre algébrique non nul sont transcendantes.

INSTITUT DE MATHÉMATIQUES, PARIS, FRANCE
*E-mail address*: andre@math.jussieu.fr

BIBLIOGRAPHIE

[Am]　Y. AMICE, Les nombres $p$-adiques, *Collection Sup.*, P.U.F. (1975).

[A1]　Y. ANDRE, *G-functions and Geometry*, Aspects of Math. **E13**, Vieweg, Braunschweig/Wiesbaden (1989).

[A2]　——, $G$-fonctions et transcendance, *J. reine angew. Math.* **476** (1996) 95–125.

[BB]　D. BERTRAND and F. BEUKERS, Equations différentielles linéaires et majorations de multiplicités, *Ann. Sci. École Norm. Sup.* **18** (1985), 181–192.

[Bé]　J. P. BEZIVIN, Sur les équations fonctionnelles aux $q$-différences, *Aequationes Math.* **43** (1992) 159–176.

[BéR]　J. P. BEZIVIN and P. ROBBA, A new $p$-adic method for proving irrationality and transcendence results, *Ann. of Math.* **129** (1989), 151–160.

[Bo]　E. BOMBIERI On $G$-functions, in *Recent Progress in Analytic Number Theory*, Vol. 2, (Durham, 1979), Academic Press, New York (1981), 1–67.

[Bun]　P. BUNDSCHUH, Quelques résultats arithmétiques sur les fonctions thêta de Jacobi, *Publ. Math. Paris* 6, *Problèmes Diophantiens* 1983-84, 1–15.

[C]　V. CHIRSKII, Relations algébriques dans les corps locaux (en russe), *Vestnik Moskov. Univ.*, sér. 1 (1990), 92–95.

[D]　B. DRAGOVICH, On some $p$-adic series with factorials, in: *p-adic Functional Analysis*, Lect. Notes in Pure and Applied Math. **192** (1997), 95-105, M. Dekker, New York.

[Sc]　W. SCHIKHOF, Ultrametric calculus, *Cambridge Studies in Adv. Math.* **4** (1984), Cambridge University Press, New York.

[Sh]　A. SHIDLOVSKII, Transcendental numbers, *de Gruyter Studies in Math.* **12** (1989); traduction de *Transtsendentnye chisla*, Nauka (1987).

[S]　C. L. SIEGEL, Über einige Anwendungen Diophantischer Approximationen, Abh. Preuss. Akad. Wiss., Phys.-Math. Kl. (1929–30), no. 1, 1–70; reproduit dans: *Gesammelte Abhandlungen* B.1, Springer (1966).

[T]　L. TSCHAKALOFF, Arithmetische Eigenschaften der unendlichen Reihe $\sum_{v=0}^{\infty} a^{-v(v-1)/2} x^v$ I: *Math. Ann.* **80** (1921) 62–74; II: *Math. Ann.* **84** (1921) 100–114.

[W]　R. WALLISSER, Rationale Appoximation des $q$-Analogons der Exponential-funktion und Irrationalitätsaussagen für diese Funktion, *Arch. Math.* **44** (1985) 59–64.

[WW]　E. WHITTAKER and G. WATSON, *A Course of Modern Analysis*, 5[th] ed., 1980, Cambridge Univ. Press, Cambridge (1[st] ed. 1902).